\input amstex
\documentstyle{ amsppt}

\magnification 1150
\NoBlackBoxes
\bigskip \noindent
\centerline {\bf $\Cal M_{15}$ IS RATIONALLY CONNECTED} \bigskip
\centerline {\bf Andrea Bruno, Alessandro Verra}
\bigskip \noindent
\bf 1. Introduction. \rm \par \noindent Let $D$ be a smooth, irreducible complex projective
curve,
it is certainly an  unexpected property that $D$ moves in a linear system
$$
\mid D \mid
$$
on a smooth irreducible surface $S$ which is not birational to $D \times \bold P^1$. For a
curve $D$ of genus $g$ with general moduli such a property is indeed equivalent to the
existence
of a rational curve $R$ such that
$$
[D] \in R \subset \Cal M_g,
$$
where $\Cal M_g$ is the moduli space of $D$ and $[D]$ denotes the moduli point of $D$. Due to
the
fundamental theorem of Eisenbud, Harris and Mumford on the Kodaira dimension of $\Cal M_g$,
there
is no rational curve through a general point of $\Cal M_g$ if $g \geq 23$. Therefore we have
not
to expect the above property for a given curve $D$ of genus $g \geq 23$. \par \noindent
Of course the existence of a rational curve through a general point of $\Cal M_g$ just means
that
$\Cal M_g$ is a uniruled variety. The uniruledness of $\Cal M_g$ is somehow expected for $g
\leq
23$, more precisely it is implied by two famous conjectures: the conjecture that every variety
of
negative Kodaira dimension is uniruled and the so called slope conjecture on effective divisors
of
$\Cal M_g$. The latter one implies that $\Cal M_g$ has negative Kodaira dimension for $g
\leq 22$, (cfr. [HM] and [FP]). \par \noindent
In spite of being expected, the uniruledness of $\Cal M_g$ for $g \leq 22$ persists to be an
open
problem if $16 \leq g \leq 22$. Let us briefly recall some history and some known
facts about this matter. For $g \leq 14$ one knows more: $\Cal M_g$ is not only uniruled but
also unirational. The proof of the unirationality of $\Cal M_g$ goes back to Severi for $g
\leq 10$, while the cases between $11$ and $14$ admit more recent proofs due to Sernesi ($g =
12$), Chang and Ran ($g = 11,13$) and
 Verra ($
g = 14$), ([S], [Se], [CR], [V]). \par \noindent  In
addition Chang and Ran showed that the Kodaira dimension of $\Cal M_g$ is negative for $g = 15,
16$ and asked about the uniruledness of these moduli spaces, ([CR1]). In this note we
positively
answer this question in the case of genus $15$. Actually we will show a stronger result:
\bigskip
\noindent
(1.1) \bf MAIN THEOREM \it $\Cal M_{15}$ is rationally connected. \rm \bigskip \noindent
The proof relies on a property which can be observed for a general curve $D$ of any genus
$g \leq 15$, namely that there exists an embedding
$$
D \subset S \subset \bold P^r
$$
where $S$ is a smooth, regular, canonical surface and $\mid D \mid
$ is at least 2-dimensional. For instance it is possible to show
that  a general $D$ of genus $g \leq 10$ admits such an embedding
in a quintic surface $S \subset \bold P^3$. Moreover it is shown
in [V] that, for $g = 11, 12, 14$, $S$ can be chosen among the few
other canonical surfaces which are also complete intersections,
(the case $g = 13$, though not covered, admits a completely
analogous description).
\par \noindent
Let $D$ be a general curve of genus $15$ and let $L$ be a general line bundle on $D$ having
degree
$9$ and satisfying $h^0(L) = 2$, it will be shown in section 3 that
$$
D \subset S \subset \bold P^6
$$
where $S$ is a complete intersection of four quadrics and $\Cal O_D(D) \cong L$. Birationally
speaking we can consider the moduli space
$$
\Cal W
$$
of pairs $(D,L)$, which is open in the universal Brill-Noether locus $\Cal W^r_{d,g}$ with $r = 1, d =
9, g
= 15$. Let
$$
h: \ \mid D \mid \to \Cal W
$$
be the natural map sending $A$ to the isomorphism class of $(A, \Cal O_A(A))$. It turns out
that $h$ is generically finite onto its image and that the same is true
for $f
\cdot h$, where
$$
f: \Cal W \to \Cal M_{15}
$$
is the forgetful map, (see 4.12). This implies a quite interesting property: \it through a
general
point of
$
\Cal W$ always passes a rational surface. \rm In particular it follows that both $\Cal W$ and
$
 \Cal M_{
15}$ are uniruled. To prove that $\Cal M_{15}$ is rationally connected we
consider the divisor
$$
\Delta_0 \subset \overline {\Cal M}_{15}
$$
parametrizing isomorphism classes of nodal stable curves of arithmetic genus 15. As is well
known
$\Delta_0$ is dominated by the moduli space $\Cal M_{14,2}$ of $2$-pointed curves of genus
$14$.
We show in section 3 that $\Cal M_{14,2}$ is unirational. Hence \it $\Delta_0$ is unirational
\rm
and the proof of the rational connectedness of $\Cal M_{15}$ can be sketched as follows: \par
\noindent let $x_i = [D_i]$, $i = 1,2$, be general in $\Cal M_{15}$ and let $L_i$ be a line
bundle
on $D_i$ with $deg \ L_i = 9$ and $h^0(L_i) = 2$. As above we have an embedding
$$
D_i \subset S_i \subset \bold P^6
$$
where $S_i$ is a smooth complete intersection of four quadrics and $\Cal O_{D_i}(D_i) \cong
L_i$.
In the 2-dimensional linear system $\mid D_i \mid$ we can choose a Lefschetz pencil $P_i$
containing $D_i$ and $D'_i$, where $D'_i$ is a nodal curve defining a \it general \rm point
$x'_i$
of $\Delta_0$. In particular this implies that $x'_i$ is smooth for $\overline {\Cal M}_{15}$,
$P_i$ defines an irreducible rational curve $R_i$ in
$\overline {\Cal M}_{15}$ containing $x_i$ and $x'_i$. Since $\Delta_0$ is unirational
$x'_1$ and $x'_2$ are connected by an irreducible rational curve $R'$. Therefore $x_1$, $x_2$
are
connected by a chain $R = R_1 \cup R_2 \cup R'$ of irreducible rational curves. Moreover
$U \cap R$ is connected, where $U$ is the regular locus of $\overline {\Cal M}_{15}$.  Hence
$\Cal M_{15}$ is rationally connected.
\par \noindent
We do not know at the moment whether $\Cal M_{15}$ is unirational, nor if this is plausible.
\par \noindent
In view of affording the uniruledness of $\Cal M_g$, in the cases where this is unknown, it
could
be interesting to have some indications about the possible embeddings $D
\subset S \subset \bold P^r$, where $S$ is a canonical surface and $D$ is general of genus $g
\in
[16,22]$. Finally it seems worth to address the following
 natural
\bigskip \noindent
(1.2) \bf PROBLEM \it For which values of $g$ is $\Cal M_g$ rationally connected? \rm \bigskip
\noindent
\bf 2. Curves of degree 19 and genus 15 in $\bold P^6$. \rm \par \noindent
In the Hilbert scheme $Hilb_{19,15,6}$ of all curves of degree $19$ and genus $15$ of $\bold
P^6$
we consider the open set $\Cal H$ whose points are stable, non degenerate, smoothable curves
and its open subset
$$
\Cal D = \lbrace D \in \Cal H \ / \ h^1(T_{\bold P^6} \otimes \Cal O_D) = 0 \rbrace. \tag 2.1
$$
Assume that
$$
\Cal D \neq \emptyset
$$
and consider any $D \in \Cal D$. As is well known the condition $h^1(T_{\bold P^6} \otimes
\Cal O_D) = 0$ implies that the Kodaira-Spencer map $dt_D: H^0(N_D) \to H^1(T_D)$ is
surjective.
Since $dt_D$ is the tangent map at $D$ of the natural map
$$
t: \Cal D \to \Cal M_{15}, \tag 2.2
$$
it follows that $t$ is dominant. On the other hand, for these values of the degree and of the
genus of $D$, there exists a unique irreducible component of $\Cal H$ which dominates $\Cal
M_{15}$. Therefore $\Cal D$ is an irreducible open subset of such a component.

We will show in a moment that $\Cal D$ is non empty. Previously we want to explain more of the
geometric situation: fix a general curve $D$ of genus $15$, then its Brill-Noether locus
$$
W^1_9(D) = \lbrace L \in Pic^9(D) \ / \ h^0(L) = 2 \rbrace \tag 2.3
$$
is a smooth, irreducible curve. This indeed follows from the general Brill-Noether theory
because the Brill-Noether number $\rho(d,g,r)$ is one if $(d,g,r) = (9,15,1)$. It is not
difficult
to check that, since
$D$ is general, the line bundle $\omega_D
\otimes L^{-1}$ is very ample and defines an embedding
$$
D \subset \bold P^6
$$
as a curve of degree $19$. $D$ is a general point in the irreducible component of $\Cal H$
which
dominates $\Cal M_{15}$. A priori we could have $h^1(T_{\bold P^6} \otimes \Cal O_D) = h^1(N_D)
\geq 1$ for every $D$ in such a component, but this is not the case if $\Cal D \neq \emptyset$.
To show that $\Cal D$ is non empty we consider the nodal, reducible
 curve
$$
D_o = R \cup A \tag 2.4
$$
where $R$ is a rational normal sextic curve in $\bold P^6$, $A$ is a smooth, irreducible,
non degenerate curve of degree $13$ and genus $8$ and
$$
R \cap A = \lbrace z_1, \dots, z_8 \rbrace =: Z \tag 2.5
$$
is a set of 8 distinct points on $R$. Note that these points are in general position, indeed
this happens for any set of $n+2$ distinct points lying on a smooth, irreducible rational
normal
curve of $\bold P^n$. Notice also that
$$
\Cal O_A(1) \cong \omega_A(-p)
$$
where $p$ is a point. We want now to show that $D_o \in \Cal D$. \bigskip \noindent
(2.6) \bf PROPOSITION \it $h^1(T_{\bold P^6} \otimes \Cal O_{D_o}) = 0$. \rm \par \noindent
PROOF Tensoring by $T_{\bold P^6}$ the standard exact sequence
$$
0 \to \Cal O_{D_o} \to \Cal O_R \oplus \Cal O_A \to \Cal O_Z \to 0
$$
and passing to the associated long exact sequence we obtain
$$
0 \to H^0(T_{\bold P^6} \otimes \Cal O_{D_o}) \to H^0(T_{\bold P^6} \otimes \Cal O_A) \oplus
H^0(T_{\bold P^6} \otimes \Cal O_R) \to \tag 2.7
$$
$$ \to H^0(T_{\bold P^6} \otimes \Cal O_Z) \to H^1(T_{\bold P^6}
\otimes \Cal O_{D_o}) \to H^1(T_{\bold P^6} \otimes \Cal O_A) \oplus H^1(T_{\bold P^6} \otimes
\Cal O_R)) \to 0.
$$
Since $\Cal O_R(1)$ is non special, the standard Euler sequence
$$
0 \to \Cal O_R \to \Cal O_R(1)^7 \to T_{\bold P^6} \otimes \Cal O_R \to 0
$$

implies that $h^1(T_{\bold P^6} \otimes \Cal O_R) = 0$. It is also standard that the
restriction
$$
b: H^0 (T_{\bold P^6} \otimes \Cal O_R) \to H^0(T_{\bold P^6} \otimes \Cal O_Z)
$$
is an isomorphism. Indeed the Euler sequence induces a diagram
$$
\CD
{H^0(\Cal O_R(1))^7} @>a>> {H^0(T_{\bold P^6} \otimes \Cal O_R)} @>b>> {H^0(T_{\bold P^6}
\otimes
\Cal O_Z))} \\
\endCD
$$
such that $a$ is surjective and $b \cdot a$ is the natural evaluation map. Then $b \cdot a$ is
surjective and $b$ is surjective, hence an isomorphism. From the previous arguments it then
suffices to show that $h^1(T_{\bold P^6} \otimes \Cal O_A)) = 0$.
Let $A'$ be the canonical model of $A$ and recall that $A$ is obtained from $A'$ by projection from $p$.
Consider the natural diagram
$$
\CD
0 @>>> H^0(A', \omega_{A'}(-p)) \otimes \Cal O_{A'}  @>>> H^0(A', \omega_{A'}) \otimes \Cal O_{A'}
@>>>  H^0(A', \omega_{A'}|_{p}) \otimes \Cal O_{A'} @>>> 0 \\
@. @VVV @VVV @VVV \\
0 @>>> \omega_{A'}(-p) @>>> \omega_{A'} @>>> \omega_{A'}|_{p}@>>> 0 \\
\endCD
$$
where the first two vertical arrows are given by evaluation and are both surjective.
If we dualize and twist the exact sequence that one obtains from the snake lemma we obtain the sequence

$$
0 \to \Cal O_{A'}(1) \to T_{\bold P^7} \otimes \Cal O_{A'}(-p) \to
T_{\bold P^6} \otimes \Cal O_A \to 0.
$$
We are then reduced to show that $h^1(T_{\bold P^7} \otimes \Cal O_{A'}(-p))) = 0$,
but this is a consequence of the fact that, since $A'$ is canonical, $h^1(T_{\bold P^7} \otimes \Cal O_{A'})=0$ and we can factor the natural surjection $H^0(\bold P^7,T_{\bold P^7}) \to T_{\bold P^7}|_p \to 0$ through $H^0(T_{\bold P^7} \otimes \Cal O_{A'}) \to T_{\bold P^7}|_p$.
\bigskip \noindent
The next theorem is a well known consequence of the vanishing of $h^1(T_{\bold P^6} \otimes
\Cal
O_{D_o})$, (cfr. [HH]1.1). \bigskip \noindent
(2.8) \bf THEOREM \it $D_o$ is smoothable. \rm \bigskip \noindent
The theorem implies that the natural morphism $t: \Cal D \to \Cal M_{15}$ is dominant. Let
$\Cal
W$ be the moduli space of pairs $(D,L)$ such that $D$ is a general curve of genus 15 and $L \in
W^1_9(D)$, in addition we can consider the natural morphism
$$
u: \Cal D \to \Cal W \tag 2.9
$$
sending $D \in \Cal D$ to the moduli point of $(D,L)$, where $L = \omega_D(-1)$. It is clear
from
the previous results and remarks that $u$ is also dominant. So we can summarize the main
achievements of this section as follows:
\bigskip \noindent
(2.10) \bf COROLLARY \it $\Cal D$ is non empty and dominates both $\Cal W$ and $\Cal M_{15}$.
\rm
 \bigskip \noindent
\bf 3. Embedding $D$ in a smooth (2,2,2,2) complete intersection of $\bold P^6$. \rm
\par
\noindent In this section we show that a general $D \in \Cal D$ is embedded in a smooth
complete
intersection of 4 quadrics and that $h^0(\Cal I_D(2)) = 4$, where $\Cal I_D$ is the ideal sheaf
of $D$. We will use a reducible curve $D_o = R \cup A$ of the type considered in the previous
section, proving that $D_o$ embeds in a reducible complete intersection of 4 quadrics $S_o$
 and
that
 the pair $(D_o,S_o)$ is smoothable. We start with a smooth, non degenerate sextic Del Pezzo
surface
$$
Y \subset \bold P^6. \tag 3.1
$$
The ideal sheaf $\Cal I_Y$ of $Y$ is generated by quadrics: this implies the next
property.
\bigskip \noindent
(3.2)\bf LEMMA \it The intersection scheme of four general quadrics containing $Y$ is a reduced
surface
$$
X \cup Y, \tag 3.3
$$
where $X$ is a smooth, irreducible component. Moreover the intersection scheme of $X$ and $Y$
is a smooth, irreducible curve
$$
B = X \cap Y. \tag 3.4
$$
\rm \par \noindent
PROOF Let $\sigma: P \to \bold P^6$ be the blowing up of $Y$, $E$ the
exceptional divisor of $\sigma$, $H$ the pull-back of a hyperplane by
$\sigma$. Since $\Cal I_Y$ is generated by quadrics the strict transform

of
$$
\mid \Cal I_Y(2) \mid
$$
is a base point free linear system and coincides with
$$
\mid 2H - E \mid.
$$
By Bertini theorem the intersection of 4 general elements $Q'_1, \dots, Q'_4 \in \ Ê\mid 2H - E
\mid $ is a smooth surface $X'$.
Since $(2H-E)^6 = 4$ is positive, $X'$ is connected. Furthermore we can assume that
$$
B' = Q'_1 \cap \dots Q'_4 \cap E
$$
is a smooth, connected curve and that $\sigma/B': B' \to \bold P^6$ is an embedding. Then
$\sigma/X': X' \to \bold P^6$ is an
embedding too and $X = \sigma(X')$ has the following properties: Ê(i) $X \cup Y$ is the
complete intersection of the quadrics
$Q_i$ $=$ $\sigma(Q'_i)$, $i = 1, \dots, 4$. (ii) $B = \sigma(B')$ is the intersection scheme
of $X$ and $Y$. This
completes the proof.
\bigskip \noindent
The complete intersection of four quadrics $X \cup Y$ is a canonical surface Êi.e. $\omega_{X
\cup Y} \cong \Cal O_{X \cup
Y}(1)$. This implies that $\omega_Y(B) \cong \omega_{X \cup Y} \otimes \Cal O_Y
\cong \Cal O_Y(1)$. Since
$\omega_Y \cong \Cal O_Y(-1)$ it follows that
$$
B \in \mid \Cal O_Y(2) \mid
$$
is a quadratic section of $Y$ and a canonical curve of genus 7. The surface
$X$ is described in detail in [V], in particular $X$ is obtained from the blowing up
$$
\sigma: X \to \bold P^2
$$
of 1
 1 points
 $e_1 \dots e_{11}$ in general position. Let $E_i$ be the exceptional divisor
over $e_i$, $L$ the pull-back of a line by $\sigma$, $H$ a hyperplane section of $X$ then
$$
\mid H \mid = \mid 6L - 2(E_1 + \dots + E_5) - E_6 - \dots - E_{11} \mid.
$$
It is easy to see that a general curve
$$
R \in \mid 2L - E_1 - E_2 - E_{10} - E_{11} \mid
$$
is a smooth, irreducible rational sextic curve. \bigskip \noindent
(3.5) \bf PROPOSITION \it $R$ is non degenerate. \rm \par \noindent
PROOF $R$ is non degenerate iff $\mid H - R \mid = \mid 4L - E_1 - E_2 - 2E_3 - 2E_4 - 2E_5 -
E_6 - E_7 - E_8 - E_9 \mid$
is empty. Note that $(H-R)E_{10} = (H-R)E_{11} = 0$ and that $(H - R)^2 = -2$. Then
consider the contraction $f: X \to X'$ of $E_{10}$ and $E_{11}$. Let $C' = f_*C$ with
$C \in \mid H - R \mid$, then ${C'}^2 = C^2 = -2$. As is well known there is no effective
$-2$ curve on the blowing up of $\bold P^2$ in $n \leq 9$ general points. Hence $C$ cannot
exist and $\mid H - R \mid$ is empty. \bigskip \noindent
Notice also that

$$
\omega_X(B) \cong \Cal O_X(H) \tag 3.6
$$
so that $H - B \sim -3L + E_1 + \dots + E_{11}$ is the canonical class.
\bigskip \noindent
(3.7) \bf PROPOSITION \it Let $\Cal I_{Y \cup R}$ be the ideal sheaf of $Y \cup R$, then
$$
h^0(\Cal J_{Y \cup R}(2)) = 4.
$$
\rm \par \noindent
PROOF Since $X \cup Y$ is the complete intersection of four quadrics and $X \cap Y = B$, it
suffices to show that $h^0(\Cal O_X(2H - B - R)) = 0$. Note that $2H - B - R \sim L - E_3 -
E_4 - E_5 + E_{10} + E_{11}$. Moreover $(2H - B - R)E_{10} =$ $(2H - B - R)E_{11} = -1$.
This implies that an effective $C \in \mid 2H - B - R \mid$ contains $E_{10}$ and
$E_{11}$. But then $C = C' + E_{10} + E_{11}$ with $C' \in \mid L - E_3 - E_4 - E_5 \mid$.
This is a contradiction because, since the points $e_3, e_4, e_5$ are not collinear, $\mid L -
E_3 - E_4 - E_5 \mid$ is empty. Hence $\mid 2H - B - R \mid$ is empty and $h^0(\Cal
O_X(2H-B-R)) = 0$.
\bigskip \noindent
(3.8) \bf PROPOSITION \it $\mid \Cal O_B(R) \mid$ is a base point free
  pencil
of
degree 8.
\rm \bigskip \noindent
PROOF One easily computes $RB = 8$. Note that $K_X - (R - B) \sim H - R$ so that, by Serre
duality, $h^i(\Cal O_X(R-B))$ $=$ $h^{2-i}(\Cal O_X(H-R))$. Since $R$ is irreducible and
$R(R-B) = -8$, we have $h^2(\Cal O_X(H-R)) = 0$. On the other hand proposition 3.5
implies that $h^0(\Cal O_X(H - R)) = 0$. Since $\chi(\Cal O_X(H-R)) = 0$, it follows
$h^1(\Cal O_X(H-R)) = 0$ and hence $h^i( \Cal O_X(R-B)) = 0$ for $i = 0,1,2$. Then
the statement follows considering the long exact sequence of
$$
0 \to \Cal O_X(R-B) \to \Cal O_X(R) \to \Cal O_B(R) \to 0
$$
\bigskip \noindent
Now we consider on $Y$ the linear system
$$
\mid B + N \mid
$$
where $N$ is one of the 6 lines contained in $Y$. We have $BN = 2$, $\ p_a(B+N) = 8$,
$\ deg \ B+N = 13$. Since $B$ is a quadratic section of $Y$ we have also $h^1(\Cal O_Y(B)) =
0$  and the exact sequence
$$
0 \to H^0(\Cal O_Y(B)) \to H^0(\Cal O_Y(B+N)) \to H^0(\Cal O_N(B)) \to 0.
$$
It easily follows from the sequence that $\mid B+N \mid$ is base point free and
that its general element is smooth, irreducible. As above we consider a general $R$ in
the pencil $\mid 2L - E_1 - E_2 - E_{10} - E_{11} \mid$. Since $\mid \Cal O_B(R) \mid$ is
base point free of degree 8 and $R$ is a rational normal sextic, we can assume that the
intersection scheme
$$
Z = B \cap R
$$
is smooth and supported on 8 points in general position. Moreover we can also assume that
$B \cap N$ and $Z$ do not intersect.  Let $\Cal I_{Z/Y}$ be the ideal sheaf of $Z$ in $Y$,
we have the following: \bigskip
\noindent
(3.9) \bf PROPOSITION \it The base locus of $\mid \Cal I_{Z/Y}(B+N) \mid$ is $Z$. \rm \bigskip
\noindent
PROOF Note that $B+N$ is smooth along $Z$ and that the linear system $\mid \Cal O_B(B+N-Z)
\mid$ is
base point free. Then the statement follows if the restriction $\rho: \mid B+N \mid \to \mid
\Cal O_B(B+N)
\mid$ is surjective. Since $h^1(\Cal O_Y(N)) = 0$, the surjectivty of $\rho$ follows from the
long exact
sequence of
$$
0 \to \Cal O_Y(N) \to \Cal O_Y(B+N) \to \Cal
O_B(B+N) \to 0.
$$
\bigskip \noindent
The proposition implies that a general
$$
A \ \in \ \mid \Cal I_Z(B+N) \mid
$$
is a smooth, irreducible curve of degree 13 and genus 8. Moreover the curve
$$
D_o = A \cup R,
$$
is nodal with $deg D_o = 19$, $p_a(D_o) = 15$, $Sing \ D_o = Z$. We know from section 2 that
$$
D_o \in \Cal D,
$$
where $\Cal D$ is the open subset of the Hilbert scheme of $D_o$ parametrizing nodal, non
degenerate, smoothable curves $D$ satisfying $h^1(T_{\bold P^6} \otimes \Cal O_D) = 0$.
\bigskip \noindent
(3.10) \bf PROPOSITION \it $h^0(\Cal I_{D_o}(2)) = 4$ and $h^i(\Cal I_{D_o}(2)) = 0$, $i > 0$.
\rm \par \noindent
PROOF We have $h^1(\Cal O_{D_o}(2)) = 0$ and $h^0(\Cal O_{D_o}(2)) = 24$, this can be easily
proved considering the standard exact sequence
$$
0 \to \Cal O_{D_o}(2) \to \Cal O_A(2) \oplus \Cal O_R(2) \to \Cal O_Z(2) \to 0
$$
and its associated long exact sequence. Then, by the standard exact sequence
$$
0 \to \Cal I_{D_o}(2) \to \Cal O_{\bold P^6}(2) \to \Cal O_{D_o}(2) \to 0,
$$
it follows that $h^i(\Cal I_{D_o}(2)) = 0$ for $i > 1$ and that
$h^1(\Cal I_{D_o}(2)) = 0$ iff $h^0(\Cal I_{D_o}(2)) = 4$. Finally
we observe that the natural inclusion $ \mid \Cal I_{Y \cup R}(2)
\mid \subseteq \mid \Cal I_{D_o}(2) \mid$ is an equality. Indeed
let $Q$ be a quadric containing $D_o$, then $Q \cap Y$ contains
$A$. But $A$ is linearly equivalent to $B+N$ where $B$ is a
quadratic section of $Y$ and $N$ is effective. Therefore $Q$
contains $Y$. Since $R \subset D_o$, it follows that $Q$ contains
$Y \cup R$ and hence $h^0(\Cal I_{D_o}(2)) = h^0(\Cal I_{Y \cup
R}(2))$. By proposition 3.7 $h^0(\Cal I_{Y \cup R}(2)) = 4$ and
this completes the proof. \bigskip \noindent The proposition
implies that a general $D \in \Cal D$ is contained in a unique
complete intersection $S$ of four quadrics. More precisely one can
show the following: \bigskip \noindent

(3.11) \bf THEOREM \it For a general $D \in \Cal D$ one has $D \subset S$, where $S$ is a
smooth complete intersection of 4 quadrics. Moreover it holds $h^0(\Cal I_D(2)) = 4$
for the ideal sheaf $\Cal I_D$ of $D$. \rm \par \noindent
PROOF The curve $D_o$ is contained in $X \cup Y$ which is a complete intersection of 4
quadrics.
Moreover, by 3.10, its ideal sheaf $\Cal I_{D_o}$ satisfies $h^0(\Cal I_{D_o}(2)) =
4$ and $h^i(\Cal I_{D_o}(2)) = 0$, $i > 0$. Then, by standard semicontinuity arguments, a
general $D \in \Cal D$ satisfies the same properties. In particular the base locus $S$ of $\mid
\Cal I_D(2) \mid$ is a complete intersection of 4 quadrics. It remains to show that $S$ is
smooth.
Let us consider the standard exact diagram of tangent and normal bundles,
$$
\CD
0 @>>> {T_{D_o}} @>>> {T_{\bold P^6}\otimes \Cal O_{D_o}} @>>> {N_{D_o}} @>>> {T^1_Z} @>>> 0 \\
@. @VVV @VVV @VVV @VVV \\
0 @>>> {T_{S_o} \otimes \Cal O_{D_o}} @>>> {T_{\bold P^6}\otimes \Cal O_{D_o}} @>>> {N_{S_o}
\otimes \Cal O_{D_o}} @>>> {T^1_{\Cal S_o} \otimes \Cal O_{D_o}} @>>> 0 \\
\endCD
$$
where $S_o = X \cup Y$ and $T^1_Z$, $T^1_{S_o}$ are the sheaves defined by the $T^1$-functor of
Lichtenbaum-Schlessinger. $T^1_Z$ and $T^1_{S_o}$ are respectively supported on the singular
loci
of $D_o$ and of $S_o$, it is easy to check that $T^1_Z = \Cal O_Z$ and that $T^1_{S_o}$ is
a line bundle on $B$. Furthermore we have $T^1_{S_o} \otimes \Cal O_{D_o} = \Cal O_{D_o \cap
B}$
and the vertical arrow
$$
 T^1_Z \to T^1_{S_o} \otimes \Cal O_{D_o}
$$
is the natural injection $\Cal O_Z \to \Cal O_{D_o \cap B}$. We know that $D_o$ is
smoothable, let $\delta \in H^0(N_{D_o})$ be an infinitesimal deformation which is induced by
an
effective smoothing of $D_o$. Then, as is well known, the image of $\delta$ in $T^1_Z$ is non
zero
at each $z \in Z$, (cfr. [S2]). Let
$$
\delta_1 \in H^0(N_{S_o} \otimes \Cal O_{D_o})
$$
be the image of $\delta$. Then, by the commutativity of the diagram, the image of
$\delta_1$ in $T^1_{S_o} \otimes \Cal O_{D_o}$ is non zero at each $z \in Z \subset D_o \cap
B$.
On the other hand every $x \in D_o \cap B - Z$ is smooth for $D_o$ and it is obvious that we
can choose
 $\delta$
 such that $\delta_1(x) \neq 0$. Finally we observe that the restriction map
$$
r: H^0(N_{S_o}) \to H^0(N_{S_o} \otimes \Cal O_{D_o})
$$
is surjective. To see this it suffices to tensor by $N_{S_o}$ the standard exact sequence
$$
0 \to \Cal I_{D_o/S_o} \to \Cal O_{S_o} \to \Cal O_{D_o} \to 0.
$$
Passing to the long exact sequence the surjectivity of $r$ follows if $h^1(\Cal I_{D_o/S_o}
\otimes N_{S_o}) = 0$. Now $N_{S_o} = \Cal O_{S_o}(2)^4$ because $S_o$ is the complete
intersection of 4 quadrics. So it suffices to show that $h^1(\Cal I_{D_o/S_o}(2)) = 0$.
This easily follows from the long exact sequence of
$$
o \to \Cal I_{S_o}(2) \to \Cal I_{D_o}(2) \to \Cal I_{D_o/S_o}(2) \to 0
$$
and proposition 3.10. Since $r$ is surjective, $\delta_1$ lifts to an infinitesimal deformation
$$
\sigma \in H^0(N_{S_o}).
$$
By the commutativity of the diagram the image of $\sigma$ in $H^0(T^1_{S_o})$ is not zero
at each $x \in B \cap D_o $. Finally let
$$
(S_t,D_t), \ t \in T,
$$
be an effective deformation of $(S_o,D_o)$ induced by $\sigma$. Then such a deformation smooths
$Z = Sing \ D_o$ so that $D_t$ is smooth. Moreover it smooths $D_o \cap B$ as a subset of $B =
Sing \ S_o$. This implies that $S_t$ has at most finitely many singular points and that they
are not in $D_t$ for $t \neq o$. Let $t \neq o$, it is easy to show that $(D_t,S_t)$ deforms to
a
pair $(D,S)$ such that $Sing \ S$ is empty: we leave this to the reader.
\bigskip \noindent
\bf 4 Proof of the main theorem. \rm \par \noindent
In this section we prove the main theorem of this paper i.e. that $\Cal M_{15}$ is rationally
connected. At first we need to show that the moduli space
$$
\Cal M_{14,2} \tag 4.1
$$
of 2-pointed curves of genus 14 is unirational. With this purpose we consider the Hilbert
scheme
$Hilb_{14,8,6}$ of curves in $\bold P^6$ having degree 14 and arithmetic genus 8. It is known
that a non degenerate, linearly normal, smooth, irreducible curve
$$
C \subset \bold P^6 \tag 4.2
$$
of degree 14 and genus 8 is projectively normal and generated by
 quadrics, ([V]). Let
$$
\Cal C \subset Hilb_{14,8,6}
$$
be the open subset parametrizing all curves with the above properties, then (see [V]Ê):
\bigskip
\noindent
(4.3) \bf PROPOSITION \it $\Cal C$ is irreducible and unirational. \rm \bigskip \noindent
It is standard to construct a projective bundle $h: P \to \Cal C$ such that the fibre
of $P$ at $C$ is
$$
P_C = \mid \Cal I_C(2) \mid,
$$
$\Cal I_C$ being the ideal sheaf of $C$. Since $C$ is projectively normal, $dim \ P_C = 6$. We
have: \bigskip \noindent
(4.4) \bf PROPOSITION \it Let $V \subset P_C$ be a general 4-dimensional linear system of a
general $C \in \Cal C$, then the base locus of $V$ is
$$
C \cup D
$$
where $D$ is a smooth, irreducible curve of genus 14 and degree 18. Conversely let $D \subset
\bold P^6$ be a curve of genus  14 and degree 18 with general moduli, then $D$ is in the base
locus of $V$ for some pair $(C,V)$.
\rm
\par
\noindent PROOF See [V]. \bigskip \noindent
Using the previous results we can easily construct a rational dominant map
$$
\phi: \Cal C \times \bold P^6 \times \bold P^6 \to \Cal M_{14,2}. \tag 4.5
$$
Indeed let $(C,x,y) \in \Cal C \times \bold P^6 \times \bold P^6$ be a sufficiently general
element, then we can assume that $x,y$ are not in $C$ and moreover that the linear system
$$
V
$$
of all quadrics containing $C \cup \lbrace x, y \rbrace$ is 4-dimensional. In view of the
previous theorem we can also assume that the base locus of $V$ is
$$
C \cup D
$$
where $D$ is a general curve of genus 14. Then we define $\phi$ by setting
$$
\phi(C,x,y) = [D,x,y] \tag 4.6
$$
where $[D,x,y]$ denotes the moduli of the 2-pointed curve $(D,x,y)$. It is clear, by
proposition
4.4, that $\phi$ is dominant. Since $\Cal C$ is unirational, it follows that \bigskip \noindent

(4.7) \bf THEOREM \it $\Cal M_{14,2}$ is unirational. \rm \bigskip \noindent
Now we consider the divisor
$$
\Delta_0 \subset \overline {\Cal M}_{15} \tag 4.8
$$
parametrizing stable singular curves of arithmetic genus 15. As is well known there exists
a natural rational
  map of
degree two
$$
\psi: \Cal M_{14,2} \to \Delta_0 \tag 4.9
$$
sending $[D,x,y]$ to the moduli point of the stable curve obtained from $D$ by glueing $x$ to
$y$.
In particular we have: \bigskip \noindent
(4.10) \bf COROLLARY \it $\Delta_0$ is unirational. \rm \bigskip \noindent
Let $D \in \Cal D$ be a general smooth curve of degree 19 and genus 15, we know from section 3
that
$$
D \subset S
$$
where $S$ is a smooth complete intersection of four quadrics. Since $D$ is non degenerate, the
sheaf $\Cal O_D(1)$ is special and $\omega_D(-1)$ is an element of the Brill-Noether locus
$W^1_9(D)$. We want to point out that $\mid \omega_D(-1) \mid$ is a base-point-free pencil and
that $D$ is linearly normal. This follows because a general $D \in \Cal D$ has general moduli.
Then $D$ has no $g^2_9$ nor a $g^1_k$ with $k \leq 8$ and hence $\mid \omega_D(-1) \mid$ is a
base-point-free $g^1_9$. Moreover, by Riemann Roch, $dim \ \mid \Cal O_D(1) \mid = 6$,
therefore
$D$ is linearly normal. Since $S$ is a canonical surface we have $\omega_D(-1) \cong \Cal
O_D(D)$.
Then, from the standard, exact sequence
$$
0 \to \Cal O_S \to \Cal O_S(D) \to \Cal O_D(D) \to 0
$$
it follows that
$$
dim \ \mid D \mid = 2. \tag 4.11
$$
Notice also that $\mid D \mid$ is base point free, because $\mid \Cal O_D(D) \mid$ is base
point
free. We expect that a general pencil in $\mid D \mid$ is a Lefschetz pencil and that a
singular
element of $\mid D \mid$ defines a general point of $\Delta_o$. We will see that this is
actually
true.
\bigskip \noindent
(4.12) \bf LEMMA \it Let $D$ be as above, then a general singular element of $\mid D \mid$ is
an irreducible curve with exactly one ordinary node and no other singularity.
\rm \par \noindent
PROOF Let $f: S \to \bold P^2$ be the covering of degree 9 defined by $\mid D \mid$. Since
$D$ is general, the linear series $\mid \Cal O_D(D) \mid$ has simple ramification. Hence the
branch curve $B$ of $f$ is reduced. This implies that a general tangent line to $B$ intersects
$B$ transversally except for the tangency point. Hence a general singular element of $\mid D \mid$
is integral with exactly one ordinary node.

\bigskip \noindent
Let us recall once more that a general smooth element $D \in \Cal D$ has the following
properties: \it
\bigskip \noindent
(i) $D \subset S$, where $S$ is a smooth complete intersection of
4 quadrics and
$h^0(\Cal I_D(2)) = 4$.  \par \noindent
(ii) $\mid \omega_D(-1) \mid$ is a base-point-free pencil and $D$ is
linearly normal. \par \noindent
(iii) The Petri map
$$
\mu_D: H^0(\omega_D(-1)) \otimes H^0(\Cal O_D(1)) \to H^0(\omega_D)
$$
is injective.
\par \noindent
(iv) The family of the singular, irreducible curves
$$
\Gamma \in \mid D \mid
$$
having an ordinary node as a unique singularity is
non empty and hence 1-dimensional. \rm \bigskip \noindent (i) holds by theorem 3.10 and (iv) by
the previous lemma, while (ii) has been just remarked above. The injectivity of $\mu_D$
follows from Gieseker-Petri theorem because $D$ has general moduli, notice also that
the Brill-Noether locus $W^1_9(D)$ is a smooth, irreducible curve. The family of all curves
$\Gamma$ as in (iv) will be denoted as
$$
\Cal D_o.
$$
(4.13) \bf LEMMA \it For each $\Gamma \in \Cal D_o$ (ii) holds for
$\Gamma$ and the Petri map $\mu_{\Gamma}$ is injective.
\rm \bigskip \noindent
PROOF Recall that $S$ is a canonical surface, so that $\omega_C(-1)$
$\cong$ $\Cal O_C(C)$ for any curve $C \subset S$. A general $D \in \mid
\Gamma \mid$ satisfies conditions (ii) and (iii), in particular $\mid
\Cal O_D(D) \mid$ is a base-point-free pencil. Then it follows from the
standard long exact sequence
$$
0 \to H^0(\Cal O_S) \to H^0(\Cal O_S(D)) \to H^0(\Cal O_D(D)) \to 0
$$
that $\mid D \mid$ is 2-dimensional and base-point-free. Replacing
$\Cal O_D$ by $\Cal O_{\Gamma}$ in the above sequence, it follows that
$\mid
\omega_{\Gamma}(-1) \mid$ is a base-point-free pencil. Let $H$ be a
hyperplane section of $S$, by Serre duality $h^1(\Cal O_S(H))$
$=$ $h^1(\Cal O_S) = 0$. Hence we have the  long exact
sequence
$$
0 \to H^0(\Cal O_S(H-D)) \to
 H^0(\Cal
 O_S(H)) \to H^0(\Cal O_D(H)) \to
H^1(\Cal O_S(H-D)) \to 0.
$$
Since a general $D \in \mid \Gamma \mid$ is non degenerate and linearly
normal, it follows
$h^0(\Cal O_S(H-D))$ $=$ $h^1(\Cal O_S(H-D))$ $ = 0$. Then, replacing
as above $\Cal O_D$ by $\Cal O_{\Gamma}$, we deduce that $\Gamma$ is non
degenerate and linearly normal. It remains to show that the Petri map
$$
\mu_{\Gamma}; H^0(\Cal O_{\Gamma}(D)) \otimes H^0(\Cal O_{\Gamma}(H)) \to
H^0(\Cal O_{\Gamma}(H+D))
$$
is injective. By the base-point-free pencil trick we have $Ker \
\mu_{\Gamma} \cong H^0(\Cal O_{\Gamma}(H-D))$, cfr. [ACGH] p.126. For the
same trick we have
$dim \ Ker \ \mu_D = h^0(\Cal O_D(H-D))$ for each $D \in \mid \Gamma
\mid$. Since $\mu_D$ is injective for a general
$D$, the map $D \to h^0(\Cal O_D(H-D))$ is generically zero on
$\mid \Gamma \mid$. Since $h^1(\Cal O_S(H-D) = 0$, the standard exact
sequence
$$
0 \to \Cal O_S(H-2D) \to \Cal O_S(H-D) \to \Cal O_D(H-D) \to 0
$$
yelds the long exact sequence
$$
0 \to H^0(\Cal O_D(H-D)) \to H^1(\Cal O_S(H-2D)) \to 0 \to \dots
$$
Then we have $h^0(\Cal
O_{\Gamma}(H-D))$
$=$
$h^1(\Cal O_S(H-2D)) = 0$ and $\mu_{\Gamma}$ is injective.
\bigskip \noindent
Let $\overline {\Cal D}$ be the closure of $\Cal D$ in the
Hilbert scheme, then $\overline {\Cal D}$ contains the dense open set
$$
U
$$
parametrizing those integral curves $D \in \Cal D$ which satisfy the
previous conditions (i), (ii), (iii), (iv) and have at most one ordinary
node as their only singularity. Note that
$$
\Cal D_o \subset U
$$
as a divisor. Indeed (i) implies that $U$ is ruled by the
family of surfaces $U_D =: U \cap \mid D \mid$, where $D \in U$. Moreover
$U_D
\cap \Cal D_0$ is pure of dimension 1. We
consider the natural map
$$
f: U \to \overline {\Cal M}_{15}.
$$
Since $U$ is open in $\Cal D$ the map $f$ is dominant, we want to show
something more: \bigskip \noindent
(4.14) \bf PROPOSITION \it The map $f/\Cal D_0: \Cal D_0 \to \Delta_0$ is
dominant. \rm \bigskip \noindent
PROOF It suffices to show that $f: U \to \overline
  {\Cal M
}_{
15}$ has
fibres of constant dimension. Then, since $\Cal D_0$ is a divisor in $U$
and $f(\Cal D_0)$ is contained in the irreducible divisor $\Delta_0$, the
statement follows from a count of dimensions. For any $D \in U$ let us consider the fibre
$F_D = f^{-1}f(D)$ and the natural morphism
$$
h: F_D \to W^1_9(D) \subset Pic^9(D)
$$
sending $D' \in F_D$ to the point $\omega_{D'}(-1)$ of the Brill-Noether
locus $W^1_9(D)$. We know that $D'$ is linearly normal, therefore $h^{-1}h(D')$ is just the
family of all curves projectively equivalent to $D'$. In particular $dim \ h^{-1}h(D')$ $=$
$dim \ PGL(7)$. On the other hand the Petri map
$$
\mu_{D'}: H^0(\omega_{D'}(-1)) \otimes H^0(\Cal O_{D'}(1)
\to H^0(\omega_{D'})
$$
is injective i.e. its corank is one. This implies that $W^1_9(D)$ is
smooth of dimension one at its point $\omega_{D'}(-1)$. Let $L$ be in
a small neighborhood $N$ of $\omega_{D'}(-1)$ in $W^1_9(D)$. Then,
by standard semicontinuity arguments, we can assume that $\mid L \mid$ is
a base-point-free pencil and that $\omega_D \otimes L^{-1}$ defines an
embedding of $D$ in $\bold P^6$ as a linearly normal curve $D"$. By
the same semicontinuity arguments we can also assume that $D" \in U$
so that $h(D") = L$. But then $N \subset h(F_D)$ and it follows that
each irreducible component of $h(F_D)$ is a curve. This implies that
$dim \ F_D =$ $dim \ PGL(7) + 1 = 49$, for each fibre $F_D$. \bigskip \noindent
(4.15) \it PROOF OF THE MAIN THEOREM \rm We can finally conclude this note by proving that
$\Cal
M_{15}$ is rationally connected. Fix two general points $[D_1]$ and $[D_2]$ in $\Cal M_{15}$,
then fix a general $L_i$ in the Brill-Noether locus $W^1_9(D_i)$, $i = 1,2$. Consider the
embedding
$$
D_i \subset \bold P^6
$$
defined by the line bundle $\omega_{D_i} \otimes L_i^{-1}$. $D_i$ is a general element of $\Cal
D$, in particular
$$
D_i \subset S_i
$$
where $S_i$ is a smooth complete intersection of 4 quadrics. Take a general pencil
$P_i \subset \mid D_i \mid$ containing $D_i$ and consider one element $D_i^0 \in P_i$ which
is
 a singular curve. Then $[D_i^0]$ is general in $\Delta_0$ and hence smooth for $\overline
{\Cal M}_{15}$. Since $\Delta_0$ is unirational, there exists an irreducible rational curve $R$
containing $[D_1]^0$ and $[D_2^0]$. Let $R_i$ be the image of $P_i$ in $\overline {\Cal
M}_{15}$
and let $U \subset \overline {\Cal M}_{15}$ be the open set of regular points. Then
$$
U \cap (R_1 \cup R_2 \cup R)
$$
is a connnected chain of rational curves joining $[D_1]$ to $[D_2]$. This implies the rational
connectedness of $\Cal M_{15}$.
\bigskip \noindent
(4.16) \bf REMARK \rm The unirationality of $\Cal W$ and $\Cal M_{15}$ would follow if there
exists a unirational variety $V \subset \Cal D$ which intersect a general $\mid D \mid$ in
finitely many points. It is not clear to us whether such a $V$ exists.

\bigskip \noindent
\bf 5. References \rm \par \noindent
[ACGH] E. Arbarello, M. Cornalba, P. Griffiths, J. Harris {\it
Geometry of Algebraic Curves I} Springer-Verlag, Berlin (1984), 1-386 \par \noindent
[CR1] M.C. Chang, Z. Ran  {\it Unirationality of the moduli space of curves of genus 11, 13
(and 12) ,} Invent. Math. {\bf 76 } (1984) 41-54,
\par \noindent
[FP] G. Farkas, M. Popa {\it Effective divisors on $M_g$ and a counterexample to the Slope
Conjecture } preprint (2002)
\par \noindent
[HH] R. Hartshorne, A. Hirschowiz
{\it Smoothing Algebraic Space Curves,} in Algebraic Geometry, Sitges 1983 (E. Casas-Alvero,
G.E. Welters, S.
Xambo-Descamps eds.)  L.N.M. {\bf 1124 } (1985), 98-131 \par
\noindent
[HM1] J. Harris, I. Morrison
{\it Moduli of Curves ,} Springer-Verlag Berlin (1991) \par \noindent
[HM2]
{\it Slopes of effective divisors on the moduli space of curves ,}  Invent. Math. {\bf 99}
(1990), 321-335 \par
\noindent
[MM] S. Mori, S. Mukai
{\it The uniruledness of the moduli space of curves of genus 11 ,} in Algebraic Geometry,
Proceedings Tokio/Kyoto 1982
L.N.M {\bf 1016 } (M.Raynaud, T. Shioda eds.) (1983), 334-353\par
\noindent
[S1] E. Sernesi  {\it L'unirazionalit\'a della variet\'a dei moduli delle
  curve di genere 12,}
Ann. Sc. Norm. Sup. Pisa {\bf 8} (1981), 405-439
\par \noindent
[S2] E. Sernesi  {\it On the existence of certain families of curves,} Invent. Math. {\bf 75 }
(1984), 25-57
\par \noindent
 [Se] F. Severi  {\it Vorlesungen ueber Algebraische Geometrie ,}  (E. Loeffler uebersetzung),
Teubner, Leipzig {\bf } (1921),\par \noindent
[ST] F.O. Schreyer, F. Tonoli  {\it Needles in a haystack: special varieties via
small fields,} in Mathematical computations with Macaulay
2,  (D. Eisenbud, D. Grayson, M. Stillman, B. Sturmfels eds.),  Springer-Verlag, Berlin (2002)
\par \noindent
[V] A. Verra {\it The unirationality of the moduli space of curves of genus $g \leq 14$}
preprint
(2004)

\end